\documentclass[11pt]{article}%
\usepackage{amsfonts}
\usepackage{mathrsfs}
\usepackage{amssymb}
\usepackage{graphicx}
\usepackage{fancyhdr}%
\usepackage{amsmath}%

\providecommand{\U}[1]{\protect\rule{.1in}{.1in}}
\renewcommand{\d}{\textrm{d}}

\newtheorem{theorem}{Theorem}[section]
\newtheorem{definition}[theorem]{Definition}
\newtheorem{lemma}[theorem]{Lemma}
\newtheorem{proposition}[theorem]{Proposition}
\begin{document}

\date{}
\title{Reflected Backward Stochastic Differential Equations with Continuous
Coefficient and $L^{2}$-Barriers}
\author{Shaolin Ji, Zhen Wu, Li Zhou\thanks{Corresponding author. \textit{E-mail
address}: zhouli@mail.sdu.edu.cn. This work is supported by the
natural Science Foundation of China (10671112), the National Basic
Research program of China (973 program, No. 2007CB814901 and No.
2007CB814904). }\\{\small {School of Mathematics, Shandong
University, Jinan 250100, PRC} }} \maketitle

\begin{abstract}
In this paper we study reflected backward stochastic differential equations
with a continuous, linear growth coefficient and two barriers which belong to
$L^{2}$. We prove that there exists at least by penalization method.
\vspace{2mm}

Keywords:\quad Backward stochastic differential equation; reflected barrier;
penalization method \vspace{1mm}

\end{abstract}

\pagestyle{fancy} \fancyhead{}
\fancyhead[CE]{Shaolin Ji , Zhen Wu , LiZhou} \fancyhead[C]{Shaolin
Ji, Zhen Wu, Li Zhou}

\vspace{2mm}

\section{Introduction}

\ \ \ \ Since Pardoux and Peng \cite{PP} introduced nonlinear backward
stochastic differential equations (BSDEs for short) with Lipschitz
coefficient, there follows many results in this topic. Lepeltier and San
Martin \cite{LS97} studied BSDEs with continuous coefficient, they proved that
in this case there exists at least one but not necessarily unique solution.
Lin and Peng. \cite{LP} got g-supersolution for BSDEs with continuous drift
coefficient. El Karoui, Kapoudjian, Pardoux, Peng, and Quenez \cite{EKPPQ}
considered reflected backward stochastic differential equations (RBSDEs for
short) for the first time, that is to say the solution should be above or
below some given process. They proved that if the coefficient is Lipschitz and
the lower barrier is continuous\thinspace, then there exists a unique
solution. And then, Lepeltier and San Martin \cite{LS2000} studied BSDEs with
continuous coefficient and two continuous barriers. In Hamad\`{e}ne \cite{H},
he studied the case of a right-continuous with left limits barrier (R.C.L.L.
for short). Recently, Lepeltier and Xu \cite{LX} gave the results of BSDEs
with Lipschitz coefficient and R.C.L.L. barriers, and then in Peng and Xu
\cite{PX} with $L^{2}$-barriers.

In this paper, we work on BSDEs with continuous coefficient and two $L^{2}%
$-barriers. We apply the result in Lepeltier and San Martin \cite{LS97}, which
showed that for a continuous function $f$, there exists a sequence of
Lipschitz function $f_{m}$ that converges to $f$ as $m\rightarrow\infty$, to
deal with the continuous coefficient. The penalization method is employed to
tackle the $L^{2}$-barriers. Our proof is also based on the monotonic limit
theorem in Peng \cite{P}.

This paper is organized as follows: in section 2, we formulate the problem for
the solutions of RBSDEs with two $L^{2}$-barriers. In section 3, some
prelilinary results are given which will be used in the proof. Then in the
last section, we give the proof of existence of solution for RBSDEs with two
$L^{2}$-barriers.

\section{Formulation of the Problem}

\ \ \ \ On a given complete probability space ($\Omega$,$\mathscr{F}$%
,P),\ \{$B_{t}$,$0\leq t\leq T$\}\ is the $d$-dimensional standard Brownain
motion,\ \{$\mathscr{F}_{t}$,$0\leq t\leq T$\} is the augmentation of the
natural filtration generated by the Brownain motion.

We introduce the following spaces:

\begin{itemize}
\item $L^{2}=\{\xi:\Omega\rightarrow\mathbb{R}^{d},\ \mathscr{F}_{T}%
\text{-measuable random variable with\ E}[|\xi|^{2}]<\infty\}$;

\item $L^{2}_{\mathscr{F}}=\{\varphi:\Omega\times[0,t]\rightarrow
\mathbb{R}^{d},\ \mathscr{F}_{t}\text{-measuable process\ with\ E}[\int
_{0}^{t}|\varphi_{t}|^{2}\text{d}t]<\infty\}$;

\item $S^{2}_{\mathscr{F}}=\{\varphi\in L^{2}_{\mathscr{F}}:
\text{progressively measurable R.C.L.L. process\ with\ E}[\sup_{0\leq t\leq
T}|\varphi_{t}|^{2}]<\infty\}$.
\end{itemize}

First of all we give the following assumptions:

\textbf{Assumption 1.} The terminal value $\xi$ is in $L^{2}$.

\textbf{Assumption 2.} The function $f:[0,T]\times\Omega\times\mathbb{R}%
\times\mathbb{R}^{d}\rightarrow f(t,w,y,z)$,\ for any $(t,w)\in[0,T]\times
\Omega,\ f(t,w,y,z)$ is continuous on $\mathbb{R}\times\mathbb{R}^{d}%
$,\ P-almost surely.\ And there exists a constant $K$,\ such that for any
$(t,y,z)\in[0,T]\times\mathbb{R}\times\mathbb{R}^{d}$,
\[
|f(t,w,y,z)|\leq K(1+|y|+|z|)\ \ \ \text{P-a.s.}
\]

\textbf{Assumption 3.} The barriers $L,U\in L_{\mathscr{F}}^{2}$ satisfy:
\[
\text{E}[\text{ess} \sup_{0\leq t\leq T}(L_{t}^{+})^{2}]<\infty,\ \text{E}%
[\text{ess} \sup_{0\leq t\leq T}(U_{t}^{+})^{2}]<\infty,\
\]
\[
L_{T}\leq\xi\leq U_{T}\ \ \text{a.s.}, \ \ \ \ \ L_{t}\leq U_{t}\ \ \text{for
all }t\in[0,T].
\]

\textbf{Assumption 4.} There exists a process
\begin{equation}
X_{t}^{0}=X_{0}^{0}+A^{0}_{t}-K_{t}^{0}+\int_{0}^{t}Z_{s}^{0}\text{d}%
B_{s}\ \ \ \ \ 0 \leq t\leq T\,,
\end{equation}
with $Z^{0}\in L_{\mathscr{F}}^{2},A^{0},K^{0}\in S_{\mathscr{F}}^{2}$,\ and
increasing with $A_{0}^{0}=K_{0}^{0}=0$, such that $L_{t}\leq X_{t}^{0}\leq
U_{t}\ \ \ $a.e.\,a.s.

We introduce the definition of the solution for RBSDE with two barriers $L$,
$U$:

\begin{definition}
A quadruple $(Y,Z,A,K)\in S_{\mathscr{F}}^{2}\times L_{\mathscr{F}}^{2} \times
S_{\mathscr{F}}^{2}\times S_{\mathscr{F}}^{2}$ is called a solution for RBSDE
with the lower barrier $L\in L_{\mathscr{F}}^{2}$, the upper barrier $U\in
L_{\mathscr{F}}^{2}$, the terminal condition $\xi\in L^{2}$ and the
coefficient $f$ if it satisfies:

\begin{enumerate}
\item $A,K$ are increasing.

\item $(Y,Z,A,K)$ satisfies the following BSDE
\begin{equation}
Y_{t}=\xi+\int_{t}^{T}f(s,Y_{s},Z_{s})\text{d}s+A_{T}-A_{t}-K_{T}+K_{t}%
-\int_{t}^{T}Z_{s}\text{d}B_{s},\ \ \ 0\leq t\leq T.
\end{equation}

\item $L_{t}\leq Y_{t}\leq U_{t}$,\ \ \ \ a.e.\ \ a.s.

\item Generalized Skorohod condition:

for each $L^{*},U^{*}\in S_{\mathscr{F}}^{2}$ such that $L_{t}\leq L_{t}%
^{*}\leq Y_{t} \leq U_{t}^{*}\leq U_{t}$\ \ \ \ a.e.\ \ a.s.,
\begin{equation}
\int_{0}^{T}(Y_{s-}-L_{s-}^{*})\text{d}A_{s} =\int_{0}^{T}(U_{s-}^{*}%
-Y_{s-})\text{d}K_{s}=0.\ \ \ \ \ \ \ \ \
\end{equation}

\end{enumerate}
\end{definition}

In this paper, our main result is the following Theorem 2.2 which will be
proved in section 4.

\begin{theorem}
\label{thoerem2} Under Assumptions (1),(2),(3),(4),\ there exists at least one
solution $(Y,Z,A,K)$ for RBSDEs with two $L^{2}$-barriers.
\end{theorem}

\section{Some Preliminary Results}

\label{basic result}

\ \ \ \ In this section,\ we introduce some preliminary definitions and
results that will be used later.\ We first introduce $g$-supersolution which
is very important for the prove of the existence theorem:

\begin{definition}
(See Peng \cite{P},\,El Karoui et al. \cite{EPQ})We call a triple $(Y,Z,A)\in
S_{\mathscr{F}}^{2} \times L_{\mathscr{F}}^{2}\times S_{\mathscr{F}}^{2}$ a
g-supersolution if $A$ is an increasing process in $S_{\mathscr{F}}^{2}$ and
the triple satisfies:
\begin{equation}
Y_{t}=Y_{T}+\int_{t}^{T} g(s,Y_{s},Z_{s})\d s+A_{T}-A_{t}-\int_{t}^{T} Z_{s}\d
B_{s}\ \ \ \ \ t\in[0,T].
\end{equation}

\end{definition}

For a continuous function with linear growth,\,we have the following lemma:

\begin{lemma}
(See Lepeltier and San Martin \cite{LS97})\ let $f:\mathbb{R}^{p}%
\rightarrow\mathbb{R}, p\in\mathbb{N},$ be a continuous function with linear
growth,\ that is to say $\forall x\in\mathbb{R}^{p},\,|f(x)|\leq K(1+|x|)$.
Define $f_{m}(x)=\inf_{y\in\mathbb{Q}^{p}}(f(y)+m|x-y|)$, then for
$m>K$,\,$f_{m}:\mathbb{R}^{p}\rightarrow\mathbb{R}$ satisfies:

\begin{enumerate}
\item Linear growth:\ $\forall x\in\mathbb{R}^{p},\,|f_{m}(x)|\leq K(1+|x|)$;

\item Monotonicity:\ $\forall x\in\mathbb{R}^{p},\,f_{m}(x)\,\uparrow f(x)$;

\item Lipschitz condition:\ $\forall x\in\mathbb{R}^{p},\,|f_{m}%
(x)-f_{m}(y)|\leq m|x-y|$;

\item Strong convergence:\ if $x_{m}\rightarrow x$,\,then $f_{m}%
(x_{m})\rightarrow f(x)$\,.
\end{enumerate}
\end{lemma}

The following generalized Monotonic Limit Theorem of BSDEs is proved in Peng
and Xu \cite{PX}.

Consider the following sequence of It\^{o}'s process:
\begin{equation}
y_{t}^{i}=y_{0}^{i}+\int_{0}^{t}g_{s}^{i}s-A_{t}^{i}+K_{t}^{i}+\int_{0}%
^{t}z_{s}^{i}B_{s},\ \ \ i=1,2,\cdots.
\end{equation}
here for each $i$,\thinspace the process $g^{i}\in L_{\mathscr{F}}^{2}%
,\,A^{i},K^{i}\in S_{\mathscr{F}}^{2}$ are given,\thinspace and $\{A^{i}%
,K^{i}\}_{i=1}^{\infty}$ satisfies
\begin{align*}
& (i)\ \ A^{i}\mbox{ is continuous and increasing such that }A_{0}%
^{i}=0\mbox{ and }\text{E}(A_{T}^{i})^{2}<\infty.\\
& (ii)\ \ K^{i}\mbox{ is increasing and }K_{0}^{i}=0.\\
& (iii)\ \ K_{t}^{i}-K_{s}^{i}\geq K_{t}^{j}-K_{s}^{j}\ \ \ \forall0\leq s\leq
t\leq T\ \ a.s.,\ \ \forall i\geq j.\\
& (iv)\ \ \mbox{for each }t\in\lbrack0,T],\,K_{t}^{j}\uparrow K_{t}\mbox{
with }\text{E}[K_{T}^{2}]<\infty.\\
&
\end{align*}
Furthermore,\thinspace we assume that
\begin{align*}
& (v)\ \ \{g^{i},z^{i}\}_{i=1}^{\infty}\mbox{ converges weakly to }(g^{0}%
,z)\mbox{ in }L_{\mathscr{F}}^{2}.\\
& (vi)\ \ \{y_{t}^{i}\}_{i=1}^{\infty}\mbox{ converges increasingly to }(y_{t}%
)\mbox{ with }\text{E}[\sup_{0\leq t\leq T}|y_{t}|^{2}]<\infty.
\end{align*}

\begin{theorem}
\label{monotonic2} Let the above assumptions hold,\thinspace we have the limit
of $\{y_{t}^{i}\}_{i=1}^{\infty}$ $(y_{t})$ has a form $y_{t}=y_{0}+\int
_{0}^{t}g_{s}^{0}s-A_{t}+K_{t}+\int_{0}^{t}z_{s}B_{s}$,\thinspace where $A$
and $K$ are increasing processes in $S_{\mathscr{F}}^{2}$.\ For each
$t\in\lbrack0,T]$\thinspace,\thinspace$A_{t}$(resp.$K_{t}$) is the
weak(resp.strong) limit of $\{A_{t}^{i}\}_{i=1}^{\infty}$(resp.$\{K_{t}%
^{i}\}_{i=1}^{\infty}$). Furthermore for any $p\in\lbrack1,2),\,\{z_{t}%
^{i}\}_{i=1}^{\infty}$ converges strongly to $z_{t}$ in $L_{\mathscr{F}}^{p}$.
\end{theorem}

\section{Proof of the Main Result}

\ \ \ \ In this section we prove Theorem \ref{thoerem2}, i.e. the existence
for the solution of RBSDEs with two $L^{2}$-barriers. Firstly, we consider,
for any integer $m$,\thinspace\ the following RBSDEs with a upper barrier
$U$:
\begin{equation}
Y_{t}^{m}=\xi+\int_{t}^{T}f_{m}(s,Y_{s}^{m},Z_{s}^{m})\text{d}s+m\int_{t}%
^{T}(L_{s}-Y_{s}^{m})^{+}\text{d}s-K_{T}^{m}+K_{t}^{m}-\int_{t}^{T}Z_{s}%
^{m}\text{d}B_{s}.
\end{equation}
Since the coefficient are Lipschitz, according to Peng and Xu \cite{PX} these
equations have unique solutions $(Y^{m},Z^{m},K^{m}),\forall m\in\mathbb{N}$.

Then for any $n,m\geq1$, we consider the following classical BSDEs:
\begin{equation}
Y_{t}^{n,m}=\xi+\int_{t}^{T}f_{m}(s,Y_{s}^{n,m},Z_{s}^{n,m})\text{d}s-\int
_{t}^{T}Z_{s}^{n,m}\text{d}B_{s}+m\int_{t}^{T}(L_{s}-Y_{s}^{n,m})^{+}%
s-n\int_{t}^{T}(Y_{s}^{n,m}-U_{s})^{+}\text{d}s.
\end{equation}
since $g_{n,m}(t,y,z)=f_{m}(t,y,z)+m(L_{t}-y)^{+}-n(y-U_{t})^{+}$ are
Lipschitz in $(y,z)$, uniformly in $(t,w)$, the equations have unique
solutions $(Y^{n,m},Z^{n,m})$. And by comparison theorem , we have that for
fixed $n,\,Y^{n,m}$ is increasing in $m$.

Set $A_{t}^{n,m}=m\int_{0}^{t}(L_{s}-Y_{s}^{n,m})^{+}\text{d}s,\ K_{t}%
^{n,m}=n\int_{0}^{t}(Y_{s}^{n,m}-U_{s})^{+}\text{d}s,\ $ we have the following proposition:

\begin{proposition}
\label{2} There exists a constant $C$ independent on $n,m$ such that
\begin{equation}
E[\sup_{0\leq t\leq T}(Y_{t}^{n,m})^{2}]+E[\int_{0}^{T}|Z_{t}^{n,m}%
|^{2}\text{d}s] +E[(A_{T}^{n,m})^{2}]+E[(K_{T}^{n,m})^{2}]\leq C.
\end{equation}

\end{proposition}

To prove this result,\ we need the following two lemmas.

Consider the following equation:
\begin{equation}
\label{9}Y_{t}^{m}=\xi+\int_{t}^{T}f_{m}(s,Y_{s}^{m},Z_{s}^{m})\text{d}s
+m\int_{t}^{T}(L_{s}-Y_{s}^{m})^{+}\text{d}s-\int_{t}^{T}Z_{s}^{m}%
\text{d}B_{s}.
\end{equation}
this is a sequence of classical BSDE,\ there exists unique solutions
$(Y^{m},Z^{m})$, for all $m \in\mathbb{N}$.

\begin{lemma}
\label{bounded} For equation (\ref{9}), we have that there exists a constant
$C$ independent of $m$ such that
\begin{equation}
E[\sup_{0\leq t\leq T}(Y_{t}^{m})^{2}]+E[\int_{0}^{T}|Z_{s}^{m}|^{2}\text{d}s]
+E[(A_{T}^{m})^{2}]\leq C.
\end{equation}
where $A_{t}^{m}:=m\int_{0}^{t}(L_{s}-Y_{s}^{m})^{+}\text{d}s$.
\end{lemma}

Apply the It\^{o}'s formula on $(Y_{t}^{m})^{2}$, the conclusion can be
deduced owe to the Gronwall's lemma and B-D-G inequality.

We can easily get a similarly result as Lemma 5.1 in Peng and Xu \cite{PX}:

\begin{lemma}
\label{lemma4} There exists a quadruple $(Y^{*},Z^{*},A^{*},K^{*})\in
S_{\mathscr{F}}^{2}\times L_{\mathscr{F}}^{2} \times S_{\mathscr{F}}^{2}\times
S_{\mathscr{F}}^{2}$ such that
\begin{equation}
\label{1}Y_{t}^{*}=\xi+\int_{t}^{T}f(s,Y_{s}^{*},Z_{s}^{*})\text{d}s+A_{T}%
^{*}-A_{t}^{*} -(K_{T}^{*}-K_{t}^{*})-\int_{t}^{T}Z_{s}^{*}\text{d}B_{s}.
\end{equation}
where $A^{*},\,K^{*}$ are both increasing, and $L_{t}\leq Y_{t}^{*}\leq U_{t}
$\ \ a.e.,a.s.
\end{lemma}

\textbf{Proof of Proposition \ref{2}:} \ Let $(Y^{+},Z^{+})$ and $(Y^{-}%
,Z^{-})$ be the solution of the following two BSDEs:
\begin{equation}
Y_{t}^{+}=\xi+\int_{t}^{T}f_{m}(s,Y_{s}^{+},Z_{s}^{+})s+A_{T}^{\ast}%
-A_{t}^{\ast}+m\int_{t}^{T}(L_{s}-Y_{s}^{+})^{+}s-\int_{t}^{T}Z_{s}^{+}B_{s}.
\end{equation}%
\begin{equation}
Y_{t}^{-}=\xi+\int_{t}^{T}f_{m}(s,Y_{s}^{-},Z_{s}^{-})s-(K_{T}^{\ast}%
-K_{t}^{\ast})-n\int_{t}^{T}(Y_{s}^{-}-U_{s})^{+}s-\int_{t}^{T}Z_{s}^{-}B_{s}.
\end{equation}
where $(A^{\ast},K^{\ast})$ is given as in Lemma \ref{lemma4}. From the
comparison theorem of the standard BSDE, we have:
\[
Y_{t}^{-}\leq Y_{t}^{n,m}\leq Y_{t}^{+}\ \ \ \forall t\in\lbrack
0,T]\ \ \text{a.s.}
\]
Review Lemma \ref{bounded}, obviously we can prove the same result if we
replace $f_{m}$ by $f$,\ or replace $m\int_{t}^{T}(L_{s}-Y_{s}^{m}%
)^{+}\text{d}s$ by $-n\int_{t}^{T}(Y_{s}^{m}-U_{s})^{+}\text{d}s$, so we
have:
\[
\text{E}[\sup_{0\leq t\leq T}(Y_{t}^{+})^{2}]+\text{E}[\sup_{0\leq t\leq
T}(Y_{t}^{-})^{2}]\leq C,
\]
then
\[
\text{E}[\sup_{0\leq t\leq T}(Y_{t}^{n,m})^{2}]\leq\max\{\text{E}[\sup_{0\leq
t\leq T}(Y_{t}^{+})^{2}],\,\text{E}[\sup_{0\leq t\leq T}(Y_{t}^{-})^{2}]\}\leq
C.
\]
For $A_{T}^{n,m}$,\thinspace\ we consider the following BSDE:
\begin{equation}
\tilde{Y}_{t}^{m}=\xi+\int_{t}^{T}f_{m}(s,\tilde{Y}_{s}^{m},\tilde{Z}_{s}%
^{m})s-(K_{T}^{\ast}-K_{t}^{\ast})+m\int_{t}^{T}(L_{s}-\tilde{Y}_{s}^{m}%
)^{+}s-\int_{t}^{T}\tilde{Z}_{s}^{m}B_{s}.\label{5}%
\end{equation}
We know $Y^{\ast}$ satisfies $L_{t}\leq Y_{t}^{\ast}\leq U_{t}$ from Lemma
4.3,\thinspace\ thus we can add the zero term $m\int_{t}^{T}(L_{s}-Y_{s}%
^{\ast})^{+}s$ to the right side of (\ref{1}).\thinspace\ Since $A_{t}^{\ast
}\geq0$,\thinspace\ from the comparison theorem,\thinspace\ it follows that
$Y_{t}^{\ast}\geq\tilde{Y}_{t}^{m}$,\thinspace\ thus $U_{t}\geq\tilde{Y}%
_{t}^{m}$,\thinspace\ then $-m\int_{t}^{T}(\tilde{Y}_{s}^{\ast}-U_{s})^{+}s$
is zero and so can be added to the right side of (\ref{5}).\thinspace\ Again
from the comparison theorem, we derive $\tilde{Y}_{t}^{m}\leq Y_{t}^{n,m}%
$,\thinspace\ and so:
\[
0\leq A_{t}^{n,m}\leq\tilde{A}_{t}^{m}:=m\int_{0}^{t}(L_{s}-\tilde{Y}_{s}%
^{m})^{+}s,
\]
then following the same process as Lemma \ref{bounded},\thinspace\ we have
$\text{E}(A_{T}^{n,m})^{2}\leq\text{E}(\tilde{A}_{T}^{m})^{2}\leq C$.

Now we consider the BSDE
\begin{equation}
\tilde{Y}_{t}^{n}=\xi+\int_{t}^{T}f(s,\tilde{Y}_{s}^{n},\tilde{Z}_{s}%
^{n})s+A_{T}^{\ast}-A_{t}^{\ast}-n\int_{t}^{T}(\tilde{Y}_{s}^{n}-U_{s}%
)^{+}s-\int_{t}^{T}\tilde{Z}_{s}^{n}B_{s}.
\end{equation}
Similarly, we can get $E(K_{T}^{n,m})^{2}\leq C$.

Apply It\^{o}'s formula to $(Y_{t}^{n,m})^{2}$, we have:
\begin{align*}
& \text{E}|Y_{t}^{n,m}|^{2}+\text{E}\int_{t}^{T}|Z_{s}^{n,m}|^{2}s\\
& \leq C(1+\text{E}\int_{t}^{T}|Y_{s}^{n,m}|^{2}\text{d}s)+\alpha\text{E}%
\int_{t}^{T}|Z_{s}^{n,m}|s+\beta\text{E}[ess\sup_{0\leq t\leq T}(L_{t}%
^{+})^{2}]\\
& +\gamma\text{E}[ess\sup_{0\leq t\leq T}(U_{t}^{-})^{2}]+\frac{1}{\beta
}\text{E}(A_{T}^{n,m})^{2}+\frac{1}{\gamma}\text{E}(K_{T}^{n,m})^{2}.
\end{align*}
choose $\alpha=\frac{1}{3}$, we get $\text{E}\int_{0}^{T}|Z_{s}^{n,m}%
|^{2}s\leq C$. The proof of Proposition 4.1 is completed.

\hfill$\Box$

To prove the Theorem \ref{thoerem2}, we let $n$ tend to $\infty$, then
\[
\left\{
\begin{array}
[c]{ll}%
Y^{n,m}\rightarrow Y^{m} & \mbox{ \ \ \ in }L_{\mathscr{F}}^{2}.\\
n\int_{0}^{T}(Y_{s}^{n,m}-U_{s})^{+}s\rightarrow K_{T}^{m} & \mbox{ \ \ \ in
}L^{2}.\\
Z^{n,m}\rightarrow Z^{m} & \mbox{ \ \ \ in }S_{\mathscr{F}}^{2}.
\end{array}
\right.
\]
where $(Y^{m},Z^{m},K^{m})$ is the unique solution of the following RBSDE:
\begin{equation}
Y_{t}^{m}=\xi+\int_{t}^{T}f_{m}(s,Y_{s}^{m},Z_{s}^{m})s-(K_{T}^{m}-K_{t}%
^{m})-m\int_{t}^{T}(L_{s}-Y_{s}^{m})^{+}s-\int_{t}^{T}Z_{s}^{m}B_{s}.\label{6}%
\end{equation}
we know that $Y^{m}\leq U$ a.e. a.s..

And we also have the following lemma:

\begin{lemma}
\label{lemma5} There exists a constant $C$ independent on $m$,\ such that
\begin{equation}
\text{E}[\sup_{0\leq t\leq T}(Y_{t}^{m})^{2}]+\text{E}[\int_{0}^{T}|Z_{t}%
^{m}|^{2}\text{d}s] +\text{E}[(A_{T}^{m})^{2}]+\text{E}[(K_{T}^{m})^{2}]\leq
C.
\end{equation}

where $A_{t}^{m}=m\int_{0}^{t}(L_{s}-Y_{s}^{m})^{+}\d s$.
\end{lemma}

From the comparison theorem, $Y^{m}$ is increasing in $m$, so there exists a
process $Y$ such that $Y^{m}\uparrow Y$,\ and from Fatou's lemma
$\text{E}[\sup_{0\leq t\leq T}Y_{t}^{2}]\leq C$.

By the dominated convergence theorem it follows that
\[
\text{E}\int_{0}^{T}(Y_{t}-Y_{t}^{m})^{2}\d t\rightarrow
0,\ \ \ \ \ \ \ \ \ \ \mbox{as}\ n\rightarrow\infty.
\]

We have already get the conclusion that $(Y^{m},Z^{m})$ is the solution of
(\ref{6}).\ Rewrite (\ref{6}) in a forward version:
\begin{equation}
Y_{t}^{m}=Y_{0}^{m}+\int_{0}^{t}f_{m}(s,Y_{s}^{m},Z_{s}^{m})s-A_{t}^{m}%
+K_{t}^{m}-\int_{0}^{t}Z_{s}^{m}B_{s}.
\end{equation}
Set $g_{t}^{m}=-f_{m}(s,Y_{s}^{m},Z_{t}^{m})$, with Lemma \ref{lemma5}, we
derive that all assumptions of Theorem \ref{monotonic2} are satisfied.\ It
follows that its limit $Y$ is in $S_{\mathscr{F}}^{2}$ and has the form
\begin{equation}
Y_{t}=\xi+\int_{t}^{T}g_{s}^{0}s+A_{T}-A_{t}-(K_{T}-K_{t})-\int_{t}^{T}%
Z_{s}B_{s}.
\end{equation}
where $(g^{0},Z,A)$ is the weak limit of $\{g(\cdot,Y^{m},Z^{m}),Z^{m}%
,A^{m}\}_{i=1}^{\infty}$ in $L_{\mathscr{F}}^{2}$,\ K is the strong limit of
$\{K_{t}^{m}\}_{i=1}^{\infty}$ in $L_{\mathscr{F}}^{2}$, $A$ and $K$ are
increasing processes in $S_{\mathscr{F}}^{2}$.\ Furthermore,\ for any
$p\in\lbrack1,2)$,\ we have $\lim\limits_{m\rightarrow\infty}E\int_{0}%
^{T}|Z_{s}^{m}-Z_{s}|^{p}s=0$.\ In Lemma 3.2,\ we showed that the sequence of
Lipschitz function $f_{m}$ converges strongly to the continuous function
$f$,\ so we get $f_{m}(\cdot,Y^{m},Z^{m})\rightarrow f(\cdot,Y,Z)$ because of
the strong convergence of $Y^{m}$ to $Y$ and the weak convergence of $Z^{m}$
to $Z$,\ and then:
\begin{equation}
Y_{t}=\xi+\int_{t}^{T}f(s,Y_{s},Z_{s})s+A_{T}-A_{t}-(K_{T}-K_{t})-\int_{t}%
^{T}Z_{s}B_{s}.\label{7}%
\end{equation}

The only problem left is to verify the generalized Skorohod
condition\thinspace.\thinspace For the upper barrier $U$,\ it is easily to
prove that for any $U^{\ast}\in S_{\mathscr{F}}^{2}$ and $Y\leq U^{\ast}\leq
U$, if we consider large enough $m$, then $Y^{m}\leq U^{\ast}\leq U$. For the
solution of the RBSDE (\ref{6}), we get $\int_{0}^{T}(U_{t-}^{\ast}-Y_{t-}%
^{m})K_{t}^{m}=0$ from the generalized Skorokhod condition. So we get
$\int_{0}^{T}(U_{t-}^{\ast}-Y_{t-})K_{t}^{m}=0$, since $0\leq U_{t-}^{\ast
}-Y_{t-}\leq U_{t-}^{\ast}-Y_{t-}^{m}$.\ Furthermore $K_{T}^{m}\uparrow K_{T}%
$,\ so $0\leq\int_{0}^{T}(U_{t-}^{\ast}-Y_{t-})(K_{t}-K_{t}^{m})\leq
(K_{T}-K_{T}^{m})\max_{t\in\lbrack0,T)}(U_{t-}^{\ast}-Y_{t-})\rightarrow
0$.\ The Skorohod condition for the upper barrier $U$ is obtained.

At last we prove the Skorohod condition holds for the lower barrier $L$.
Consider the following BSDE:
\begin{equation}
\tilde{Y}_{t}^{m}=\xi+\int_{t}^{T}f_{m}(s,\tilde{Y}_{s}^{m},\tilde{Z}_{s}%
^{m})s+m\int_{t}^{T}(L_{s}-\tilde{Y}_{s}^{m})^{+}s-(K_{T}-K_{t})-\int_{t}%
^{T}\tilde{Z}_{s}^{m}B_{s}.
\end{equation}
We denote $\bar{Y}^{m}:=\tilde{Y}^{m}-K$ and rewrite the BSDE:
\begin{equation}
\bar{Y}_{t}^{m}=\xi-K_{T}+\int_{t}^{T}f_{m}^{K}(s,\bar{Y}_{s}^{m},\tilde
{Z}_{s}^{m})s+m\int_{t}^{T}(L_{s}-K_{s}-\bar{Y}_{s}^{m})^{+}s-\int_{t}%
^{T}\tilde{Z}_{s}^{m}B_{s}.
\end{equation}
where $f_{m}^{K}(t,y,z):=f_{m}(t,y+K,z)$.

If we consider a BSDE with coefficient $f^{K}$ and lower barrier $L^{K}$,
where $f^{K}(t,y,z)=f(t,y+K,z),\ L^{K}=L-K$, then the BSDE above is the
penalized equation of this problem, we know that it has the unique solution
$(\bar{Y}^{m},\tilde{Z}^{m},\tilde{A}^{m})$. When $m\rightarrow\infty$, we get
the limit:
\begin{equation}
\bar{Y}_{t}=\xi-K_{T}+\int_{t}^{T}f^{K}(s,\bar{Y}_{s},\tilde{Z}_{s}%
)s+\tilde{A}_{T}-\tilde{A}_{t}-\int_{t}^{T}\tilde{Z}_{s}B_{s}.\label{8}%
\end{equation}
here $\tilde{A}_{t}$ is the $L_{\mathscr{F}}^{2}$ weak limit of $\tilde{A}%
_{t}^{m}=m\int_{0}^{t}(L_{s}-\tilde{Y}_{s}^{m})^{+}s=m\int_{0}^{t}(L_{s}%
-K_{s}-\bar{Y}_{s}^{m})^{+}s$. Suppose $\ddot{Y}$ is another $f^{K}%
$-supersolution with decomposition $(\ddot{Z},\ddot{A})$,\ which satisfies
(\ref{8}) and $\ddot{Y}_{t}\geq L_{t}-K_{t}$.\ By comparison theorem,\ we have
$\bar{Y}_{t}^{m}\leq\ddot{Y}_{t}$,\ so $\bar{Y}_{t}\leq\ddot{Y}_{t}$.\ That is
to say $\bar{Y}$ is the smallest $f^{K}$-supersolution with $\bar{Y}_{T}%
=\xi-K_{T}$ that dominates $L-K$,\thinspace and from the comparison theorem we
have $Y_{t}^{m}\geq\tilde{Y}_{t}^{m}$, so we get:
\[
\tilde{A}_{t}^{m}-\tilde{A}_{s}^{m}=m\int_{s}^{t}(L_{r}-\tilde{Y}_{r}^{m}%
)^{+}dr=A_{t}^{m}-A_{s}^{m}\ \ \ 0\leq s\leq t\leq T.
\]
thus $\tilde{A}_{t}-\tilde{A}_{t}\geq A_{t}-A_{s}$.\ From (\ref{7}) we know
$Y-K$ is a $f^{K}$-supersolution,\ compare this with (\ref{8}),\ we have
$Y-K\leq\bar{Y}$.\ thus $Y-K=\bar{Y}$ is the smallest $f^{K}$-supersolution
with terminal condition $\xi-K_{T}$that dominates $L-K$.\ With the help of the
following Proposition 4.6, we can get that for each $L^{\ast}\in
S_{\mathscr{F}}^{2}$ such that $Y\geq L^{\ast}\geq L$,\thinspace we have
$Y-K\geq L^{\ast}-K\geq L-K$, then:
\[
\int_{0}^{T}(Y_{t-}-L_{t-}^{\ast})A_{t}=\int_{0}^{T}((Y_{t-}-K_{t-}%
)-(L_{t-}^{\ast}-K_{t-}))A_{t}=0.
\]
The proof of the existence for solution of RBSDEs is completed.\thinspace

\hfill$\Box$

\ \ At last, we prove that if $Y$ is the smallest $f$-supersolution that
dominates $L$, then $Y$ satisfies the Skorohod condition, which was used in
above proof. \

According to Peng and Xu \cite{PX},\ we have the following proposition:

\begin{proposition}
\label{Y-Y:} Given $Y\in S_{\mathscr{F}}^{2}$,\ $Y_{T}=\xi\in L^{2},\ L\in
L_{\mathscr{F}}^{2}$,\ the following two items are equivalent:

i) $Y$ is the smallest $g$-supersolution that dominates $L$.

ii) For any $L^{*}\in S_{\mathscr{F}}^{2},\,Y_{t}\geq L_{t}^{*}\geq L_{t}%
$,\,a.e.,a.s.,\ $Y$ is the smallest $g$-supersolution that dominates $L^{*}$.
\end{proposition}

Now we consider the following condition:\ $L\in S_{\mathscr{F}}^{2}$ is a
given process,\ $f_{0}(t)\equiv0$,\ $\hat Y\in S_{\mathscr{F}}^{2}$ is a
$f_{0}$-supersolution that dominates $L$ with terminal condition $\xi$,\ i.e.
\begin{equation}
\hat Y_{t}=\xi+A_{T}-A_{t}-\int_{t}^{T} Z_{s}\d B_{s},\ \ \ \hat Y_{t}\geq
L_{t},\ \ \ \forall t\in[0,T]\ \ \text{a.s.}%
\end{equation}
where $(Z,A)$ is the corresponding composition of $\hat Y$.\ From Peng and Xu
\cite{PX},\ we know that if $\hat Y$ is the smallest $f_{0}$-supersolution
that dominates $L$ with terminal condition $\xi$,\ then for each stopping time
$\tau\leq T$,\ we have $\hat Y_{\tau-}=\hat Y_{\tau}\vee L_{\tau-}$.\ Then we
have:
\begin{equation}
\label{0GSC}\sum_{0\leq t\leq T}(\hat Y_{t-}-L_{t-})(A_{t}-A_{t-}%
)=0\ \ \ \text{a.s.}%
\end{equation}

\begin{proposition}
\label{Y=Y:} We claim that the following two items are equivalent:

i) $Y$ is the smallest $f$-supersolution that dominates $L$ with terminal
condition $\xi$.

ii) $\hat{Y}$ is the smallest $f_{0}$-supersolution that dominates $\hat{L}$
with terminal condition $\hat{\xi}$,\ where for each $t\in[0,T]$:
\[
\hat{f}(t):=f(t,Y_{t},Z_{t}),\ \ \ \hat{Y}_{t}:=Y_{t}+\int_{0}^{t}\hat{f}(s)\d
s,\ \ \ \hat{L}_{t} :=L_{t}+\int_{0}^{t} \hat{f}(s)\d s,\ \ \ \hat{\xi}%
:=\xi+\int_{0}^{T}\hat{f}(s)\d s.
\]

\end{proposition}

\textbf{Proof:} We consider the following penalized BSDE:
\begin{equation}
\tilde{Y}_{t}^{m}=\xi+\int_{t}^{T}\hat{f}(s)s+m\int_{t}^{T}(L_{s}-\tilde
{Y}_{s}^{m})^{+}s-\int_{t}^{T}\tilde{Z}_{s}^{m}B_{s}.
\end{equation}
Comparing it with the penalized BSDE:
\begin{equation}
\hat{Y}_{t}^{m}=\xi+\int_{t}^{T}\hat{f}(s)s+m\int_{t}^{T}(L_{s}+\int_{0}%
^{t}f(s)s-\hat{Y}_{s}^{m})^{+}s-\int_{t}^{T}\tilde{Z}_{s}^{m}B_{s}.
\end{equation}
we know that we only need to prove $\hat{Y}_{t}^{m}\rightarrow\hat{Y}$,\ then
we have $\tilde{Y}_{t}^{m}\rightarrow Y$.

Suppose $\{(\tilde{Y}^{m},\tilde{Z}^{m})\}_{m=1}^{\infty}$ converges to
$(\tilde{Y},\tilde{Z})$,\ then $\tilde{Y}$ is the smallest $\hat{f}%
$-supersolution that dominates $L$ with terminal condition $\xi$.\ Next we
prove that $(\tilde{Y},\tilde{Z})=(Y,Z)$.\ Apply It\^{o}'s formula to
$|Y_{t}^{m}-\tilde{Y}_{t}^{m}|^{2}$,\ we have:
\begin{align*}
\text{E}|Y_{t}^{m}-\tilde{Y}_{t}^{m}|^{2}+\text{E}\int_{t}^{T}|Z_{s}%
^{m}-\tilde{Z}_{s}^{m}|^{2}\d s  & =2\text{E}\int_{t}^{T}(Y_{s}^{m}-\tilde
{Y}_{s}^{m})(f_{m}(s,Y_{s}^{m},Z_{s}^{m})-\hat{f}(s))\d s\\
& +2m\text{E}\int_{t}^{T}(Y_{s}^{m}-\tilde{Y}_{s}^{m})((L_{s}-Y_{s}^{m}%
)^{+}-(L_{s}-\tilde{Y}_{s}^{m})^{+})\d s.
\end{align*}
It's easy to check that $(Y_{s}^{m}-\tilde Y_{s}^{m})((L_{s}Y_{t}^{m}%
)^{+}-(L_{s}-\tilde{Y}_{s}^{m})^{+})\leq0$. Then we have:
\begin{align*}
\text{E}|Y_{t}^{m}-\tilde{Y}_{t}^{m}|^{2}+\text{E}\int_{t}^{T}|Z_{s}%
^{m}-\tilde{Z}_{s}^{m}|^{2}\d s  & \leq2\text{E}\int_{t}^{T}(Y_{s}^{m}%
-\tilde{Y}_{s}^{m})(f_{m}(s,Y_{s}^{m},Z_{s}^{m})-\hat{f}(s))\d s\\
& \leq2\text{E}\int_{t}^{T}(|Y_{s}^{m}-Y_{s}|+|\tilde{Y}_{s}^{m}-\tilde{Y}%
_{s}|)|f_{m}(s,Y_{s}^{m},Z_{s}^{m})-\hat{f}(s)|\d s\\
& +2\text{E}\int_{t}^{T}(Y_{s}-\tilde{Y}_{s})(f_{m}(s,Y_{s}^{m},Z_{s}%
^{m})-\hat{f}(s))\d s.
\end{align*}

We know $|Y^{m}-Y|+|\tilde Y^{m}-\tilde Y|\rightarrow0$ in $L_{\mathscr{F}}%
^{2}$ and $|f_{m}(s,Y_{s}^{m},Z_{s}^{m})-\hat{f}(s)|$ is uniformly bounded in
$L_{\mathscr{F}}^{2}$. Moreover from the strong convergence of $\{Y^{m}%
\}_{m=1}^{\infty}$ to $Y$ and weak convergence of $\{Z^{m}\}_{m=1}^{\infty}$
to $Z$ in $L_{\mathscr{F}}^{2}$ , we know that $\{f_{m}(\cdot,Y^{m}%
,Z^{m})\}_{m=1}^{\infty}$ converges weakly to $\hat{f}(\cdot)$. Thus the right
side of the above inequality converges to zero, it follows that $\tilde
Y\equiv Y$ and $\tilde Z\equiv Z$.

\hfill$\Box$

Follow the Theorem 4.1 d) $\Rightarrow$ e) in Peng and Xu \cite{PX}, we can
derive directly that $\hat{Y}$ defined in Proposition \ref{Y=Y:} ii) satisfies
the following condition: for each $\hat{L}^{\ast}\in S_{\mathscr{F}}^{2}$ such
that $\hat{Y}\leq\hat{L}^{\ast}\leq\hat{L}$, a.e. a.s.,
\begin{equation}
\int_{0}^{T}(\hat{Y}_{t-}-\hat{L}_{t-}^{\ast})A_{t}=0\ \ \ \ \ \text{a.s.}%
\label{S-C_}%
\end{equation}
Thus, we get the Skorohod condition of $\hat{Y}$ which is our desired result.

\end{document}